\documentclass[12pt,reqno]{amsart}
\usepackage{amssymb}
\usepackage{amsfonts}
\usepackage{graphicx}
\usepackage{color}
\usepackage{hyperref}
\usepackage{booktabs}

\newtheorem{theorem}{Theorem}[section]
\newtheorem{proposition}[theorem]{Proposition}
\newtheorem{definition}[theorem]{Definition}

\newtheorem{corollary}[theorem]{Corollary}
\newtheorem{example}[theorem]{Example}

\def\F{\mathbb{F}}
\def\Z{\mathbb{Z}}
\def\P{\mathcal{P}}
\def\D{\mathcal{D}}
\DeclareMathOperator{\dev}{dev}
\DeclareMathOperator{\Aut}{Aut}
\DeclareMathOperator{\Atop}{Atop}
\DeclareMathOperator{\Apar}{Apar}

\begin{document}

\title{Projection cubes of symmetric designs}

\author[V.~Kr\v{c}adinac and L.~Reli\'{c}]{Vedran Kr\v{c}adinac and Lucija Reli\'{c}}

\address{Faculty of Science, University of Zagreb, Bijeni\v{c}ka cesta~$30$, HR-$10000$ Zagreb, Croatia}

\email{vedran.krcadinac@math.hr}
\email{lucija.relic@math.hr}

\thanks{This work has been supported by the Croatian Science Foundation
under the project $9752$.}



\date{April 8, 2025}

\begin{abstract}
We introduce a new type of $n$-dimensional generalization of symmetric $(v,k,\lambda)$
block designs. We prove upper bounds on the dimension~$n$ in terms of~$v$ and~$k$. We
also define the corresponding concept of $n$-dimensional difference sets, and extend some
classic families of difference sets to higher dimensions. Complete classifications
are performed for small parameters $(v,k,\lambda)$ and some interesting examples
are presented.
\end{abstract}

\maketitle

\section{Introduction}

In~\cite{dL90} and~\cite{dLH93}, a wide class of $n$-dimensional combinatorial
designs was studied, including Hadamard matrices, symmetric block designs,
(generalized) weighing matrices, and orthogonal designs. Higher-di\-men\-si\-o\-nal
versions of these designs are defined as $n$-dimensional matrices such that
all $2$-dimensional sections (submatrices obtained by fixing all but two
coordinates) satisfy the appropriate condition. E.g.\ for Hadamard matrices,
the $2$-dimensional sections are $v\times v$ matrices~$H$ with
$\{-1,1\}$-entries satisfying $H\cdot H^t = vI$.

Higher-dimensional symmetric designs of this type were recently studied in~\cite{KPT24}
under the name \emph{cubes of symmetric designs}. Throughout this paper, symmetric
designs are identified with their incidence matrices. Thus, a symmetric $(v,k,\lambda)$
design is a $v\times v$ matrix~$A$ with $\{0,1\}$-entries satisfying
$A\cdot A^t = (k-\lambda)I+\lambda J$, where $I$ is the identity
matrix and~$J$ is the all-ones matrix. For the usual definition as incidence
structures of points and blocks, and background material on symmetric designs,
we refer to~\cite{BJL99} and~\cite{IT07}.

Other types of combinatorial designs have been generalized to higher
dimensions in different ways. For example, Room cubes are $n$-di\-men\-si\-o\-nal
matrices such that every $2$-dimensional \emph{projection} is a Room
square~\cite{JHD07}. By using projections instead of sections, equivalence
of Room squares and other combinatorial objects carries over to higher
dimensions. Among other things, Room squares are equivalent to pairs
of orthogonal $1$-factorizations of the complete graph, as well as pairs
of orthogonal-symmetric latin squares. Room cubes of dimension~$n$ are
equivalent to~$n$ mutually orthogonal such objects.

In the present paper, we study an analogous generalization of symmetric
designs: a \emph{$(v,k,\lambda)$ projection $n$-cube} is an $n$-dimensional
matrix such that every $2$-di\-men\-si\-onal projection is a symmetric
$(v,k,\lambda)$ design. We find that the definition using projections
works remarkably well in this context. An important and widely studied
construction technique for symmetric designs are difference sets~\cite{JPS07}.
The definition of a difference set extends nicely to~$n$ dimensions, suitable
for constructing projection $n$-cubes. Some classic families of difference
sets obtained from finite fields can be generalized to higher dimensions.
For cubes of symmetric designs in the sense of~\cite{KPT24}, the dimension
can be arbitrarily large for given parameters $(v,k,\lambda)$. On the other
hand, the dimension of a $(v,k,\lambda)$ projection cube is bounded by~$v$
and~$k$, similarly as for Room cubes.

The layout of our paper is as follows. In Section~\ref{sec2}, we define
$(v,k,\lambda)$ projection $n$-cubes and prove some of their basic properties.
Proposition~\ref{charoa} establishes equivalence to a subclass of orthogonal
arrays $OA(vk,n,v,1)$, giving a practical representation of projection cubes.
The special cases $k=1$ and $k=2$ are completely described, and examples are
given for parameters $(7,3,1)$. Up to isomorphism, the Fano plane is the unique
$(7,3,1)$ design and can be seen as a projection cube of dimension $n=2$. For
$n=3$, there are already inequivalent $(7,3,1)$ projection cubes.
Theorem~\ref{dimbound} gives an upper bound on the dimension~$n$
in terms of~$v$. Theorem~\ref{dimbound2} gives a better bound
in terms of~$v$ and~$k$.

Section~\ref{sec3} introduces $n$-dimensional difference sets. Projection
cubes coming from difference sets are characterized as having an autotopy
group acting sharply transitively on each coordinate. An example is given, showing
that they need not have any additional symmetries. In contrast, the difference
cubes from~\cite{KPT24} are invariant under a larger autotopy group and, for
abelian groups, under any conjugation.

In Section~\ref{sec4}, we generalize the Paley difference sets, the cyclotomic
difference sets, and the twin prime power difference sets to higher dimensions.
The development of the Paley difference set in the field of order~$7$ is a
$7$-dimensional analog of the Fano plane.

Computational results about projection cubes with small parameters $(v,k,\lambda)$
are presented in Section~\ref{sec5}. We perform complete classifications
of $n$-dimensional difference sets, and determine the exact value of the highest
possible dimension $\mu(v,k,\lambda)$ for some parameters. For $(16,6,2)$,
we construct $3$-dimensional projection cubes that cannot be obtained from
difference sets. Four interesting examples are presented, having non-isomorphic
$(16,6,2)$ designs as projections.

Finally, in Section~\ref{sec6}, we give some concluding remarks and outline
future research directions.

\section{Definition and basic properties}\label{sec2}

An $n$-dimensional matrix of order~$v$ over~$\F$ is a function $$C:\{1,\ldots,v\}^n\to \F,$$
where $\{1,\ldots,v\}^n$ stands for the Cartesian $n$-ary power. For $1\le x<y\le n$, we
define the projection $\Pi_{xy}(C)$ as the $2$-dimensional matrix with $(i_x,i_y)$-entry
\begin{equation}\label{projsum}
\sum_{1\le i_1,\ldots,\widehat{i_x},\ldots,\widehat{i_y},\ldots,i_n\le v} C(i_1,\ldots,i_n).
\end{equation}
The sum is taken over all $n$-tuples $(i_1,\ldots,i_n)\in \{1,\ldots,v\}^n$
with fixed coordinates $i_x$ and $i_y$.

\begin{definition}\label{defpc}
A \emph{$(v,k,\lambda)$ projection $n$-cube} is a matrix
$$C:\{1,\ldots,v\}^n \to \F$$ with $\{0,1\}$-entries such
that all projections $\Pi_{xy}(C)$, $1\le x<y\le n$ are
symmetric $(v,k,\lambda)$ designs. The set of all such matrices
will be denoted by $\P^n(v,k,\lambda)$.
\end{definition}

We can interpret~$C$ as a characteristic function and identify it
with the subset of $n$-tuples
$$\overline{C}=\{(i_1,\ldots,i_n)\in \{1,\ldots,v\}^n \mid C(i_1,\ldots,i_n)=1\}.$$

\begin{example}\label{exC12}
Two $(7,3,1)$ projection $3$-cubes $C_1$ and $C_2$ are given by
\begin{align*}
\overline{C}_1 =\{ & (1,2,3),\, (1,4,5),\, (1,6,7),\, (2,3,1),\, (2,4,6),\, (2,7,5),\, (3,1,2),\\
 & (3,6,5),\, (3,7,4),\, (4,3,7),\, (4,5,1),\, (4,6,2),\, (5,1,4),\, (5,2,7),\\
 & (5,3,6),\, (6,2,4),\, (6,5,3),\, (6,7,1),\, (7,1,6),\, (7,4,3),\, (7,5,2)\},
\end{align*}
\begin{align*}
\overline{C}_2 =\{ & (1,2,7),\, (1,4,3),\, (1,6,5),\, (2,3,6),\, (2,4,5),\, (2,7,1),\, (3,1,4),\\
 & (3,6,2),\, (3,7,5),\, (4,3,1),\, (4,5,2),\, (4,6,7),\, (5,1,6),\, (5,2,4),\\
 & (5,3,7),\, (6,2,3),\, (6,5,1),\, (6,7,4),\, (7,1,2),\, (7,4,6),\, (7,5,3)\}.
\end{align*}
\end{example}

\begin{figure}[t]
\begin{center}
\includegraphics[width=127mm]{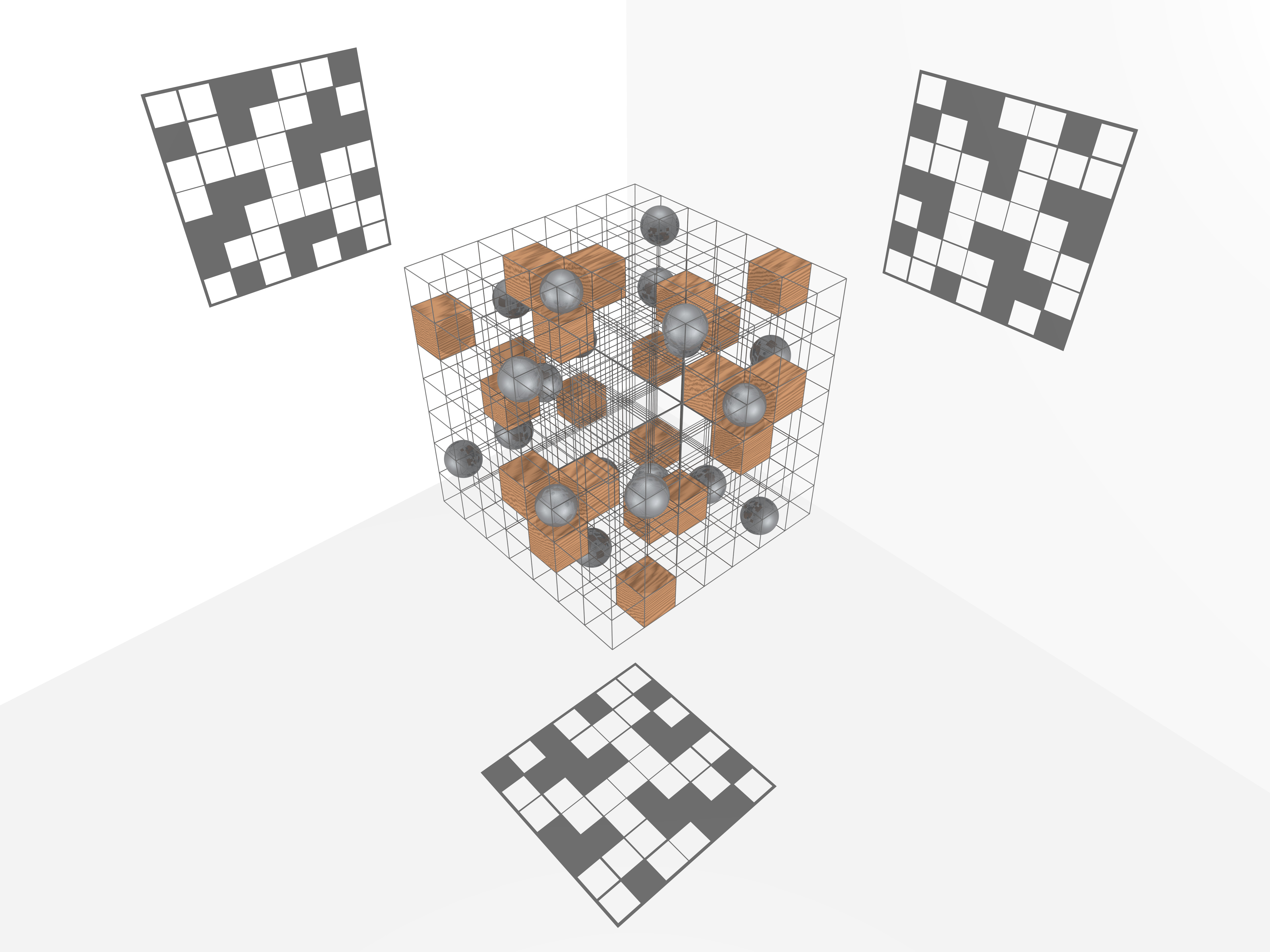}
\end{center}
\caption{Two disjoint cubes in $\P^3(7,3,1)$.}\label{fig1}
\end{figure}

The sets $\overline{C}_1$ and $\overline{C}_2$ are disjoint. The corresponding $3$-cubes
are depicted in Figure~\ref{fig1}, with incidences ($1$-entries) of~$C_1$ shown as cubes,
and of~$C_2$ as spheres. The image was rendered using the ray tracing software
POV-Ray~\cite{POVRay}. Three light sources were placed along the coordinate axes
to make the projections visible as shadows. Notice that both $3$-cubes have
identical projections: $\Pi_{ij}(C_1)=\Pi_{ij}(C_2)$ for $(i,j)\in \{(1,2),\,(1,3),\,(2,3)\}$.
If the projections are interpreted as physical shadows, the sum $C_1+C_2$ (or union
$\overline{C}_1\cup \overline{C}_2$) would also be a $(7,3,1)$ projection cube.
Moreover, any sphere in Figure~\ref{fig1} could be deleted without affecting the shadows,
giving $2^{21}$ examples of such objects. This interpretation corresponds to taking
the sum~\eqref{projsum} in the binary semifield, with $1+1=1$. Taking the sum in the
binary field~$\F_2$, with $1+1=0$, would also give many examples. In this case
switching $0\leftrightarrow 1$ in any subset of coordinates with an even number
of entries in each direction would not change the projections.

To avoid such examples with the number of incidences not determined by the
parameters, $\F$ is assumed to be a field of characteristic~$0$ in
Definition~\ref{defpc}. Then, each sum~\eqref{projsum} contains at most
one entry~$1$, because the projections $\Pi_{xy}(C)$ are $\{0,1\}$-matrices.
The incidence matrix of a $(v,k,\lambda)$ design contains $vk$ entries~$1$,
therefore the total number of incidences in~$C$ is also~$vk$. We record this as

\begin{proposition}
For each $C\in \P^n(v,k,\lambda)$, the cardinality of $\overline{C}$ is $vk$.
\end{proposition}

Projections of subsets $S\subseteq \{1,\ldots,v\}^n$ are defined as restrictions to
pairs of coordinates:
$$\Pi_{xy}(S)=\{(i_x,i_y) \mid (i_1,\ldots,i_n)\in S\}.$$
For an $n$-dimensional matrix $C:\{1,\ldots,n\}^n\to \{0,1\}$ such that
$\Pi_{xy}(C)$ is a $\{0,1\}$-matrix, $\overline{\Pi_{xy}(C)}=\Pi_{xy}(\overline{C})$ holds.
The following proposition characterizes subsets corresponding to $(v,k,\lambda)$
projection cubes.

\begin{proposition}\label{charoa}
Let $S\subseteq \{1,\ldots,v\}^n$ be a subset of cardinality $vk$. There exists
a $C\in\P^n(v,k,\lambda)$ such that $S=\overline{C}$ if and only if the following
statements are true for all $1\le x<y\le n$:
\begin{enumerate}
\item for all $i\in \{1,\ldots,v\}$, there are exactly~$k$ elements $j\in \{1,\ldots,v\}$
such that $(i,j)\in \Pi_{xy}(S)$,
\item for all $j\in \{1,\ldots,v\}$, there are exactly~$k$ elements $i\in \{1,\ldots,v\}$
such that $(i,j)\in \Pi_{xy}(S)$,
\item for all $i,i'\in \{1,\ldots,v\}$, $i \neq i'$, there are exactly~$\lambda$
elements $j\in \{1,\ldots,v\}$ such that $(i,j)\in \Pi_{xy}(S)$ and $(i',j)\in \Pi_{xy}(S)$.
\end{enumerate}
\end{proposition}

\begin{proof}
The three statements are equivalent to $\Pi_{xy}(S)$ being the set of
incidences of a symmetric $(v,k,\lambda)$ design. The assumption
$|S|=vk$ excludes examples with more incidences, like
$S=\overline{C}_1\cup \overline{C}_2$.
\end{proof}

\begin{corollary}\label{oacor}
If $C$ is a $(v,k,\lambda)$ projection $n$-cube, then $\overline{C}$ is an
orthogonal array of size~$vk$, degree~$n$, order~$v$, strength~$1$, and
index~$k$, i.e.\ an $OA(vk,n,v,1)$.
\end{corollary}

See~\cite{GC07} or~\cite{HSS99} for the definition of general orthogonal
arrays. The particular case $OA(vk,n,v,1)$ means that each element of
$\{1,\ldots,v\}$ appears exactly~$k$ times in every coordinate of~$\overline{C}$.
This follows directly from property~$(1)$ or~$(2)$ of Proposition~\ref{charoa}.

Next, we define equivalence and symmetry of projection cubes. The terminology
and definitions are the same as in~\cite{KPT24}. The cubes $C,C'\in
\P^n(v,k,\lambda)$ are \emph{isotopic} if there exist permutations
$\alpha_1,\ldots,\alpha_n\in S_v$ such that
$$C'(i_1,\ldots,i_n) = C(\alpha_1^{-1}(i_1),\ldots,\alpha_n^{-1}(i_n)),
\kern 2mm \forall (i_1,\ldots,i_n) \in \{1,\ldots,v\}^n.$$
Then, $\overline{C'}=\{(\alpha_1(i_1),\ldots,\alpha_n(i_n)) \mid (i_1,\ldots,i_n)\in
\overline{C}\}$ and the $n$-tuple of permutations $(\alpha_1,\ldots,\alpha_n)$ is
called an \emph{isotopy} between~$C$ and~$C'$, or \emph{autotopy} if $C=C'$. The set
of all autotopies of~$C$ forms a group with coordinatewise composition, called the
\emph{full autotopy group} $\Atop(C)$. \emph{Conjugation} by $\gamma\in S_n$ means
permuting the order of the coordinates:
$$C^\gamma (i_1,\ldots,i_n) = C(i_{\gamma^{-1}(1)},\ldots,i_{\gamma^{-1}(n)}),\kern 2mm
\forall (i_1,\ldots,i_n) \in \{1,\ldots,v\}^n.$$
The cubes are \emph{equivalent} or \emph{paratopic} if $C'$ is isotopic to a
conjugate of~$C$.

The two cubes of Example~\ref{exC12} are not isotopic, but
they are equivalent. The $(1,2)$-conjugate of $\overline{C}_1$ is transformed
into~$\overline{C}_2$ by applying the isotopy $\alpha_1=\alpha_2=\alpha_3=(2,3,4,5,6,7)$.
Here is an example not equivalent to~$C_1$ and~$C_2$:

\begin{example}\label{exC3}
A cube $C_3\in \P^3(7,3,1)$ is given by
\begin{align*}
\overline{C}_3 =\{ & (1,2,4),\, (1,4,6),\, (1,6,2),\, (2,3,7),\, (2,4,3),\, (2,7,4),\, (3,1,6),\\
 & (3,6,7),\, (3,7,1),\, (4,3,6),\, (4,5,3),\, (4,6,5),\, (5,1,2),\, (5,2,3),\\
 & (5,3,1),\, (6,2,7),\, (6,5,2),\, (6,7,5),\, (7,1,4),\, (7,4,5),\, (7,5,1)\}.
\end{align*}
\end{example}

Our cubes are distinguished by the \emph{full autoparatopy group} $\Apar(C)$, i.e.\
all combinations of conjugation and isotopy mapping~$C$ onto itself. Equivalent
cubes have isomorphic full autoparatopy groups, while our examples have $|\Apar(C_1)|=|\Apar(C_2)|=63$
and $|\Apar(C_3)|=42$. The GAP~\cite{GAP} package \emph{Prescribed Automorphism Groups}~\cite{PAG}
contains commands \texttt{Cube\-Filter} and \texttt{CubeAut} to decide equivalence and
compute auto\-topy and auto\-paratopy groups of $n$-cubes. Implementation details are
given in \cite[Section~2]{KPT24}, and examples of how to use these and other commands
for working with projection cubes are available in the PAG manual~\cite{PAG}.

Symmetric designs with $k=1$ are a degenerate case: $(v,1,0)$ designs exist
for all positive integers~$v$. Likewise, every orthogonal array $OA(v,n,v,1)$
fulfils the conditions of Proposition~\ref{charoa} and represents a
$(v,1,0)$ projection $n$-cube. This is a set of the form
$$\{ (i,\alpha_2(i),\ldots,\alpha_n(i)) \mid i=1,\ldots,v\}$$
for arbitrary permutations $\alpha_2,\ldots,\alpha_n\in S_v$. Hence,
$|\P^n(v,1,0)|=(v!)^{n-1}$ and all these cubes are equivalent.

\begin{table}[!b]
\caption{Number of $\P^n(3,2,1)$-cubes up to equivalence.}\label{tab1}
\begin{tabular}{cccccc}
\toprule
$n$: & 2 & 3 & 4 & 5 & 6\\
\midrule
$\#$: & 1 & 2 & 1 & 1 & 0\\
\bottomrule
\end{tabular}
\end{table}

Symmetric designs with $k=2$ exist only for $(3,2,1)$. A complete classification
is easily performed for these small parameters. The result is presented in
Table~\ref{tab1}. There are no $\P^n(3,2,1)$-cubes for $n\ge 6$, and for
$n\le 5$ the following example describes all cubes up to equivalence.

\begin{example}\label{exC45}
A $3$-cube $C_4\in \P^3(3,2,1)$ is given by
\begin{align*}
\overline{C}_4 =\{ & (1,1,1),\, (1,2,2),\, (2,1,3),\, (2,3,2),\,
(3,2,3),\, (3,3,1)\}.
\end{align*}
A $5$-cube $C_5\in \P^5(3,2,1)$ is given by
\begin{align*}
\overline{C}_5 =\{ & (1,1,1,1,1),\, (1,2,2,2,2),\, (2,1,2,3,3),\, (2,3,3,1,2),\\
 & (3,2,3,3,1),\, (3,3,1,2,3)\}.
\end{align*}
Restrictions of $\overline{C}_5$ to any~$n$ coordinates are
$\P^n(3,2,1)$-cubes of dimensions $n=2$, $3$, and $4$. These
restrictions are mutually equivalent for fixed~$n$, and for
$n=3$ they are not equivalent to~$C_4$.
\end{example}

Let $\nu(v,k,\lambda)$ be the largest integer~$n$ such that
a $(v,k,\lambda)$ projection $n$-cube exists. We have just seen
that $\nu(v,1,0)=\infty$ and $\nu(3,2,1)=5$. The next theorem
shows that $\nu(v,k,\lambda)$ is finite for all $k\ge 2$.

\begin{theorem}\label{dimbound}
If a $(v,k,\lambda)$ projection $n$-cube with $k\ge 2$ exists,
then $$n\le \frac{v(v+1)}{2}.$$
\end{theorem}

\begin{proof}
Let $\overline{C}$ be the orthogonal array representation of a cube
$C\in \P^n(v,k,\lambda)$. For any two $n$-tuples
$a,b\in \overline{C}$ and any pair of coordinates $1\le x < y \le n$,
we have $(a_x,a_y)\neq (b_x,b_y)$. This follows because both~$\overline{C}$
and the projection $\Pi_{xy}(\overline{C})$ contain $vk$ elements.

Take any $a^{(1)}\in \overline{C}$ with $a^{(1)}_1=1$. By applying
isotopy to coordinates~$2$ to~$n$, we can make $a^{(1)}=(1,1,\ldots,1)$.
Next, take any $a^{(2)}\in \overline{C}$ with $a^{(2)}_1=2$.
This $n$-tuple contains at most one entry~$1$ in coordinates~$2$
to~$n$. Assuming $a^{(2)}_x=a^{(2)}_y=1$ for some $2\le x<y\le n$
gives $(a^{(1)}_x,a^{(1)}_y) = (a^{(2)}_x,a^{(2)}_y) = (1,1)$,
a contradiction. If a $1$-entry does appear in~$a^{(2)}$, we can move
it to coordinate~$2$ by applying conjugation. Coordinates~$3$ to~$n$
contain entries different from~$1$, and by applying isotopy we
can bring $a^{(2)}$ into the form $(2,*,2,\ldots,2)$ without
affecting $a^{(1)}$. Now take $a^{(3)}\in \overline{C}$ such that
$a^{(3)}_1=3$. This element contains at most one entry~$1$ and at
most one entry~$2$ in coordinates~$3$ to~$n$. By applying conjugation
and isotopy, we can bring it into the form $a^{(3)}=(3,*,*,*,3,\ldots,3)$
without affecting~$a^{(1)}$ and~$a^{(2)}$. Continuing in the same way,
we get elements in $\overline{C}$ of the form $a^{(4)}=(4,*,*,*,*,*,*,4,\ldots,4)$,
$\ldots$, $a^{(v)}=(v,*,\ldots,*,v,\ldots,v)$. The element $a^{(v)}$
contains $1+2+3+\ldots+(v-1)=v(v-1)/2$ unknown entries~$*$.

Since $k\ge 2$, there is another element $b\in \overline{C}$, $b\neq a^{(1)}$
with $b_1=1$. Coordinates $v(v-1)/2+2$ to~$n$ of $b$ are distinct from~$1$
because of~$a^{(1)}$, and they are mutually distinct because of
$a^{(2)},\ldots,a^{(v)}$. Hence, $n\le 1+v(v-1)/2+v-1$ holds and the
inequality follows.
\end{proof}

For parameters $(3,2,1)$, Theorem~\ref{dimbound} gives the bound 
$\nu(3,2,1)\le 6$ which is not tight. An anonymous referee has 
pointed out that a better bound can be obtained from coding theory:

\begin{theorem}\label{dimbound2}
If a $(v,k,\lambda)$ projection $n$-cube with $k\ge 2$ exists,
then $$n\le \frac{vk-1}{k-1}.$$
\end{theorem}

\begin{proof}
The orthogonal array representation $\overline{C}$ of a
cube $C\in \P^n(v,k,\lambda)$ is a code of length~$n$ over an
alphabet of~$v$ symbols. We will determine its distance
distribution with respect to the Hamming metric
$$d(a,b) = | \{ x\in \{1,\ldots,n\} \mid a_x \neq b_x\} |.$$
For a fixed codeword $a\in \overline{C}$ and $i=0,\ldots,n$,
define
$$A_i=| \{b\in \overline{C} \mid d(a,b)=i\} |.$$
Then, $A_0=1$ and $\sum_{i=0}^n A_i = |\overline{C}|=vk$ clearly
hold. We claim that $d(a,b)\ge n-1$ for all $b\in \overline{C}$,
$b\neq a$. Otherwise, there would exist coordinates $1\le x<y\le n$
such that $(a_x,a_y)=(b_x,b_y)$, and $\Pi_{xy}(C)$ would not
be a matrix of zeros and ones. Therefore, $A_i=0$ for
$i=1,\ldots,n-2$. Next, $A_{n-1}=n(k-1)$ follows from
Corollary~\ref{oacor}, i.e.\ from the fact that each symbol
appears exactly~$k$ times in every coordinate. We can now compute
$A_n$ as $$A_n=vk - \sum_{i=0}^{n-1} A_i = vk-1 - n(k-1).$$
The inequality follows from $A_n \ge 0$ and $k\ge 2$.
\end{proof}

Theorem~\ref{dimbound2} gives a tight bound $\nu(3,2,1)\le 5$
in the smallest case.

\section{Higher-dimensional difference sets}\label{sec3}

Let $G$ be an additively written group of order~$v$, not necessarily
abelian. We can index projection cubes with elements of~$G$ instead
of the integers $\{1,\ldots,v\}$. Then, an $n$-cube is a function
$C:G^n\to \{0,1\}$, or the corresponding subset $\overline{C}\subseteq G^n$.
A \emph{$(v,k,\lambda)$ difference set} in~$G$ is a subset $D\subseteq G$
of size~$k$ such that every element $g\in G\setminus\{0\}$ can be written
as $g=d_1-d_2$ for exactly~$\lambda$ choices of $d_1,d_2\in D$. Its
\emph{development} $\dev D =\{ D+g \mid g\in G \}$ is the set of
blocks of a symmetric $(v,k,\lambda)$ design. The group~$G$ is
an automorphism group of the design acting sharply transitively
on the points and blocks by right translation.

An incidence matrix of $\dev D$ can be written using the Iverson
bracket: $$C:G^2\to \{0,1\},\kern 3mm C(g,h)= [g\in D+h] = [g-h\in D].$$
The expression $[P]$ takes the value~$1$ if $P$ is true and~$0$ otherwise,
see~\cite{DEK92}. The corresponding set of incident pairs is
$\overline{C}=\{(g,h)\in G^2 \mid g-h\in D\} =
\{ (g,-d+g) \mid g\in G, d\in D\}$. We can think of the difference
set as pairs $\{(0,-d) \mid d\in D\}$ and write $\overline{C}$
as a development: $\overline{C}=\{ (0,-d)+g \mid g\in G, d\in D\}$.
Using this notation, difference sets extend naturally to higher
dimensions.

\begin{definition}\label{defdifset}
An \emph{$n$-dimensional $(v,k,\lambda)$ difference set} in~$G$ is a
set of $n$-tuples $D\subseteq G^n$ of size~$k$, such that $\{d_x-d_y
\mid d\in D\}\subseteq G$ are $(v,k,\lambda)$ difference sets for
all $1\le x < y\le n$.
\end{definition}

\begin{proposition}\label{dzprop}
If $D$ is an $n$-dimensional $(v,k,\lambda)$ difference set in~$G$, then
the development
$$\dev D = \{(d_1+g,\ldots,d_n+g) \mid g\in G,\, d\in D \}$$
is the representation $\overline{C}\subseteq G^n$ of a projection cube
$C\in \P^n(v,k,\lambda)$.
\end{proposition}

\begin{proof}
First, notice that $|\dev D|=vk$. The development clearly does not contain
more than $vk$ elements. Assuming $d+g=d'+g'$ for some $g,g'\in G$ and $d,d'\in D$
implies $g=g'$ and $d=d'$, because the sets $\{d_x-d_y \mid d\in D\}$ contain~$k$
distinct elements of~$G$.

We now show that $S=\dev D$ satisfies the properties of Proposition~\ref{charoa} for
all $1\le x<y\le n$. For every $i\in G$, we count elements $j\in G$ such that
$(i,j)\in \Pi_{xy}(S)$, i.e.\ $(i,j)=(d_x+g,d_y+g)$ for some $g\in G$, $d\in D$.
This is equivalent to $j=-(d_x-d_y)+i$. Again, there are exactly~$k$ such elements
because $|\{d_x-d_y \mid d\in D\}|=k$. Thus, property~$(1)$ holds, and property~$(2)$
is checked analogously.

It remains to verify property~$(3)$ of Proposition~\ref{charoa}. For every pair of distinct
elements $i,i'\in G$, we need to show that there are exactly~$\lambda$ elements $j\in G$
such that $(i,j),(i',j)\in \Pi_{xy}(S)$. An easy calculation shows that this is equivalent
to $(d_x-d_y)-(d_x'-d_y')=i-i'$ for exactly~$\lambda$ choices $d,d'\in D$. The latter
holds because $\{d_x-d_y \mid d\in D\}$ is a $(v,k,\lambda)$ difference set.
\end{proof}

Exchanging any $n$-tuple $d\in D$ with a translate $d+g=(d_1+g,\ldots,d_n+g)$
does not affect the defining property from Definition~\ref{defdifset}, nor
does it change the development $\dev D$. In this way, an $n$-dimensional difference
set can be normalized so that its $n$-tuples have $0$ in the first coordinate.
Then the restriction of $D$ to any other coordinate is an ``ordinary'' $(v,k,\lambda)$
difference set (subset of $G$). Conversely, an ordinary difference set can be made
into a $2$-dimensional difference set by adding a coordinate containing~$0$.

\begin{example}\label{exD12}
Let $G$ be the group $\Z_7=\{0,\ldots,6\}$ with addition modulo~$7$. Then,
$D_1=\{( 0, 1, 3 )$, $(0, 2, 6)$, $( 0, 4, 5)\}$ is a $3$-dimensional
$(7,3,1)$ difference set such that $\dev D_1$ is equivalent to the
projection cubes~$\overline{C}_1$ and~$\overline{C}_2$ of Example~\ref{exC12}.
The development of $D_2=\{ (0, 1, 2 )$, $( 0, 2, 4 )$, $( 0, 4, 1 )\}$ is
equivalent to~$\overline{C}_3$ of Example~\ref{exC3}.
\end{example}

Difference sets in arbitrary groups were first considered in~\cite{RHB55}.
There it was proved that the ``left differences'' $-d_1+d_2$, $d_1,d_2\in D$
also cover $G\setminus\{0\}$ exactly~$\lambda$ times. From this, it follows
that the left translates $\{g+D\mid g\in G\}$ are the blocks of another
$(v,k,\lambda)$ design. In general, this design is isomorphic to the
dual of~$\dev D$. Of course, for abelian groups the left and right
translates coincide.

We note that for $n$-dimensional difference sets $D\subseteq G^n$
the left translates $\{ (g+d_1,\ldots,g+d_n) \mid g\in G, d\in D\}$
need \emph{not} represent projection $n$-cubes. For this, the left
differences of coordinates $\{-d_x+d_y \mid d\in D\}$ would have to
be $(v,k,\lambda)$ difference sets for all $1\le x < y\le n$.
Here is a counterexample in multiplicative notation:

\begin{example}\label{exD3}
Let $G$ be the finitely presented non-abelian group
$$\langle a,b\mid a^4=b^4=1,\, ba=ab^3\rangle.$$
This is a semidirect product $\Z_4\rtimes \Z_4$ of order~$16$. Let
$$D_3=\{(1,1),\, (1,a^2),\, (1,b^3),\, (1,a^2b),\, (1,ab^2),\, (a^2b,a^3b)\} \subseteq G^2.$$
The right differences $\{d_1 d_2^{-1} \mid d\in D_3\}=
\{1,a^2,b,a^2b^3,a^3b^2,a^3\}$ are a $(16,6,2)$ difference
set in~$G$. Hence, $D_3$ is a $2$-dimensional difference set,
and $\dev D_3$ is the set of incident pairs of a symmetric
$(16,6,2)$ design. However, $\{d_1^{-1} d_2 \mid d\in D_3\}$ is not
a difference set, because the left differences of $(1,ab^2)$
and $(a^2b,a^3b)$ coincide. Therefore, their left translates also
coincide: $\{g\cdot(1,ab^2) \mid g\in G\} = \{g\cdot(a^2b,a^3b) \mid g\in G\}$.
Thus, there are not enough incident pairs in the ``left development''
of~$D_3$ for a $(16,6,2)$ design.
\end{example}

\begin{proposition}\label{difatop}
Let $D$ be an $n$-dimensional $(v,k,\lambda)$ difference set in~$G$.
Then, the projection cube $\overline{C}=\dev D$ has an autotopy group
isomorphic to~$G$, acting sharply transitively on each coordinate.
\end{proposition}

\begin{proof}
For $g\in G$, denote right translation by $\alpha_g:G\to G$, $\alpha_g(x)=x+g$.
This is a permutation of~$G$. The set $\overline{G}=\{(\alpha_g,\ldots,\alpha_g)
\mid g\in G\}$ with coordinatewise composition is a group isomorphic to~$G$.
Acting on coordinates, this group leaves $\dev D$ invariant and is therefore
a subgroup of $\Atop(C)$. The action on each coordinate is right translation,
which is sharply transitive.
\end{proof}

Cubes of symmetric designs from~\cite{KPT24} of arbitrary dimension~$n$
can be constructed from ordinary difference sets $D\subseteq G$, see
\cite[Theorem~3.1]{KPT24}. These cubes have a much larger autotopy
group $G^{n-1}$ \cite[Theorem~3.4]{KPT24}. If~$G$ is abelian, they are
invariant under any conjugation \cite[Proposition~3.3]{KPT24}. In
contrast, projection cubes constructed from $n$-dimensional difference
sets need not have any additional symmetries apart from
Proposition~\ref{difatop}, even for abelian groups~$G$.

\begin{example}\label{exD4}
Let~$G$ be the direct product $\Z_4\times \Z_4$. A $3$-dimensional
$(16,6,2)$ difference set in~$G$ is given by
\begin{align*}
D_4 = \{ & ((0,0),(0,0),(1,0)),\, ((0,0),(1,0),(0,0)),\, ((0,0),(0,1),(2,0)),\\[1mm]
 & ((0,0),(2,0),(0,1)),\, ((0,0),(1,2),(0,3)),\, ((0,0),(2,3),(3,2))\}.
\end{align*}
A computation in PAG~\cite{PAG} shows that the full autoparatopy group
of $\dev D_4$ is of order~$16$. The differences of coordinates
$\{d_x-d_y \mid d\in D_4\}$, $(x,y)\in\{(1,2),(1,3),(2,3)\}$
are three pairwise inequivalent $(16,6,2)$ difference sets in~$G$.
This means that they cannot be mapped onto each other by automorphisms
of~$G$ and translations.
\end{example}

\begin{proposition}
Let $C\in \P^n(v,k,\lambda)$ be a projection cube with an autotopy
group~$G$ acting sharply transitively on each coordinate. Then, there
is an $n$-dimensional $(v,k,\lambda)$ difference set $D\subseteq G^n$
such that $\overline{C}$ is equivalent to $\dev D$.
\end{proposition}

\begin{proof}
Assuming that~$C$ is indexed by $\{1,\ldots,v\}$, number the group
elements $G=\{g^{(1)},\ldots,g^{(v)}\}$. Each element is an autotopy
$g^{(i)}=(\alpha_1^{(i)},\ldots,\alpha_n^{(i)})$ with permutations
$\alpha_j^{(i)}\in S_v$ as components. Since~$G$ acts sharply
transitively on the first coordinate, the elements can be numbered 
so that $\alpha_1^{(i)}(1)=i$, $\forall i\in\{1,\ldots,v\}$.
By applying isotopy to~$C$, we may also achieve $\alpha_j^{(i)}(1)=i$,
$\forall i\in\{1,\ldots,v\}$ for components $j=2,\ldots,n$. We can now
identify $i\leftrightarrow g^{(i)}$ so that $C: G^n\to \{0,1\}$
and~$G$ acts by right translation on the coordinates. This action
splits $\overline{C}\subseteq G^n$ into~$k$ orbits of size~$v$.
Take any set of orbit representatives $D\subseteq G^n$. Then, $|D|=k$
and $\dev D = \overline{C}$ obviously hold.

We claim that $D$ is an $n$-dimensional $(v,k,\lambda)$ difference set.
Take any $1\le x<y\le n$ and look at the set $\{d_x-d_y\mid d\in D\}\subseteq G$
(the operation in~$G$ is coordinatewise composition, but we write it
additively). We know that $\Pi_{xy}(\dev D)$ satisfies property~$(3)$ of
Proposition~\ref{charoa}: for all distinct $i,i'\in G$, there are
exactly $\lambda$ elements $j\in G$ such that $(i,j),(i',j)\in
\Pi_{xy}(\dev D)$. Just like in the proof of Proposition~\ref{dzprop},
this is equivalent to $(d_x-d_y)-(d'_x-d'_y)=i-i'$ for exactly~$\lambda$
choices of $d,d'\in D$. Hence, $\{d_x-d_y\mid d\in D\}$ is a $(v,k,\lambda)$
difference set.
\end{proof}

\section{Constructions from finite fields}\label{sec4}

Let $q\equiv 3 \pmod{4}$ be a prime power and $\F_q^*$ the non-zero elements
in the finite field of order~$q$. It is known that the squares in~$\F_q^*$
constitute a $(q,(q-1)/2,(q-3)/4)$ difference set in $(\F_q,+)$,
as well as the non-squares. These are the \emph{Paley difference sets}, and we
denote them by $\F_q^\square$
and~$\F_q^{\raisebox{0.7pt}{$\scriptstyle\not\mathrel{\raisebox{-0.7pt}{$\scriptstyle\square$}}$}}$.
The next theorem gives a higher-dimensional version of Paley difference sets.

\begin{theorem}\label{tmpaley}
If $q\equiv 3 \pmod{4}$ is a prime power, then there exists a $q$-dimensional
difference set with parameters $(q,(q-1)/2,(q-3)/4)$ in the additive
group of~$\F_q$.
\end{theorem}

\begin{proof}
Let $\alpha$ be a primitive element of $\F_q$. The elements of~$\F_q^\square$
are powers of $\alpha$ with even exponents, and the elements
of~$\F_q^{\raisebox{0.7pt}{$\scriptstyle\not\mathrel{\raisebox{-0.7pt}{$\scriptstyle\square$}}$}}$
with odd exponents. Consider the following subset of~$(\F_q)^q$:
\begin{equation}\label{eqpaley}
D=\{ (0,\alpha^{2i},\alpha^{2i+1},\alpha^{2i+2},\ldots,\alpha^{2i+q-2}) \mid i=0,\ldots,(q-3)/2\}.
\end{equation}
The differences of elements in the first coordinate and another coordinate~$y$
are $\{-\alpha^{2i+y} \mid i=0,\ldots,(q-3)/2\}$. Since $-1$ is not a square,
this set is $\F_q^\square$ for odd~$y$ and
$\F_q^{\raisebox{0.7pt}{$\scriptstyle\not\mathrel{\raisebox{-0.7pt}{$\scriptstyle\square$}}$}}$
for even~$y$. Similarly, the differences of two non-zero coordinates $x$, $y$ are
$\{\alpha^{2i+x}-\alpha^{2i+y} \mid i=0,\ldots,(q-3)/2\} =
\{\alpha^{2i+x}(1-\alpha^{y-x}) \mid i=0,\ldots,(q-3)/2\}$. This is
again either $\F_q^\square$
or~$\F_q^{\raisebox{0.7pt}{$\scriptstyle\not\mathrel{\raisebox{-0.7pt}{$\scriptstyle\square$}}$}}$,
depending on the parity of~$x$ and whether $1-\alpha^{y-x}$ is a square or non-square.
Since the differences of any two coordinates make a $(q,(q-1)/2,(q-3)/4)$ difference set,
$D$ is a $q$-dimensional difference set according to Definition~\ref{defdifset}.
\end{proof}

\begin{example}\label{exD5}
Substituting $q=7$ and $\alpha=3$ in~\eqref{eqpaley}, we get the following
$7$-dimensional $(7,3,1)$ difference set in $\Z_7$:
$$D_5=\{ (0, 1, 3, 2, 6, 4, 5 ),\,\, ( 0, 2, 6, 4, 5, 1, 3 ),\,\, ( 0, 4, 5, 1, 3, 2, 6 )\}.$$
Thus, the Fano plane extends at least to dimension~$7$ as a projection cube: $\nu(7,3,1)\ge 7$.
Theorem~\ref{dimbound} gives the upper bound $\nu(7,3,1)\le 28$, and Theorem~\ref{dimbound2}
gives $\nu(7,3,1)\le 10$. Notice that the $3$-dimensional difference sets~$D_1$ and~$D_2$ from
Example~\ref{exD12} are restrictions of $D_5$ to coordinates $\{1,2,3\}$ and $\{1,2,4\}$,
respectively.
\end{example}

The fourth powers in $\F_q^*$ and their cosets are $(q,(q-1)/4,(q-5)/16)$ difference sets
for $q=4t^2+1$, $t$ odd. The eighth powers and their cosets are $(q,(q-1)/8,(q-9)/64)$
difference sets for suitable~$q$ (e.g.~$q=73$), see~\cite{LB71} and~\cite{JPS07}.
These \emph{cyclotomic difference sets} can also be made into $q$-dimensional versions:
$$D=\{ (0,\alpha^{mi},\alpha^{mi+1},\alpha^{mi+2},\ldots,\alpha^{mi+q-2}) \mid i=0,\ldots,(q-1)/m-1\}$$
where $m=4$ or $m=8$. The proof is very similar as for Theorem~\ref{tmpaley}, which is the case $m=2$.

The \emph{twin prime power difference sets} are another classic family constructed from finite
fields. Suppose that both $q$ and $q+2$ are odd prime powers and let $\alpha$
and $\beta$ be primitive elements of $\F_q$ and $\F_{q+2}$, respectively. Let $G=\F_q \times \F_{q+2}$
be the direct product of additive groups and look at the following subsets:
\begin{align}
\begin{split}\label{twinds1}
E_1 &= \{(\alpha^{2i},\beta^{2j})\mid i=0,\ldots,(q-3)/2,\, j=0,\ldots,(q-1)/2\},\\[1mm]
E_2 &= \{(\alpha^{2i+1},\beta^{2j+1})\mid i=0,\ldots,(q-3)/2,\, j=0,\ldots,(q-1)/2\},\\[1mm]
E_3 &= \{ (0,0) \} \cup \{(\alpha^i,0) \mid i=0,\ldots,q-2\}.
\end{split}
\end{align}
In~\cite[Theorem~VI.8.2]{BJL99} it is proved that $D=E_1\cup E_2\cup E_3$ is a difference set
in~$G$ with Paley-Hadamard parameters $(4m-1,2m-1,m-1)$ for $m=(q+1)^2/4$. The construction can
be modified by taking
\begin{align}
\begin{split}\label{twinds2}
E'_1 &= \{(\alpha^{2i},\beta^{2j+1})\mid i=0,\ldots,(q-3)/2,\, j=0,\ldots,(q-1)/2\},\\[1mm]
E'_2 &= \{(\alpha^{2i+1},\beta^{2j})\mid i=0,\ldots,(q-3)/2,\, j=0,\ldots,(q-1)/2\}.
\end{split}
\end{align}
Then, $D'=E'_1\cup E'_2\cup E_3$ is also a difference set in~$G$, with the same parameters.
The following theorem gives a generalization to higher dimensions.

\begin{theorem}\label{tmtpp}
If $q$ and $q+2$ are odd prime powers, then there exists a $q$-dimensional
difference set in $G=\F_q \times \F_{q+2}$ with parameters $(4m-1,2m-1,m-1)$
for $m=(q+1)^2/4$.
\end{theorem}

\begin{proof}
Take the following sets of $q$-tuples from $G^q$:
\begin{align*}
E_1 &= \{(0,a^{(i,j)}_0,\ldots,a^{(i,j)}_{q-2})\mid i=0,\ldots,(q-3)/2,\, j=0,\ldots,(q-1)/2\},\\[1mm]
E_2 &= \{(0,a^{(i,j)}_1,\ldots,a^{(i,j)}_{q-1})\mid i=0,\ldots,(q-3)/2,\, j=0,\ldots,(q-1)/2\},\\[1mm]
E_3 &= \{ (0,\ldots,0) \} \cup \{(0,b^{(i)}_0,\ldots,b^{(i)}_{q-2}) \mid i=0,\ldots,q-2\}.
\end{align*}
Here, $0$ is the neutral element in~$G$ and $a^{(i,j)}_{x} = (\alpha^{2i+x},\beta^{2j+x})$,
$b^{(i)}_{x}=(\alpha^{i+x},0) \in G$. We claim that $D=E_1\cup E_2\cup E_3$ is a $q$-dimensional
difference set, i.e.\ that the differences of any pair of coordinates give a $(4m-1,2m-1,m-1)$
difference set in~$G$. By taking the first coordinate and subtracting another coordinate~$y$, we
get three sets
\begin{align*}
 & \{-(\alpha^{2i+y},\beta^{2j+y})\mid i=0,\ldots,(q-3)/2,\, j=0,\ldots,(q-1)/2\},\\[1mm]
 & \{-(\alpha^{2i+y+1},\beta^{2j+y+1}) \mid i=0,\ldots,(q-3)/2,\, j=0,\ldots,(q-1)/2\},\\[1mm]
 & \{(0,0)\}\cup\{-(\alpha^{i+y},0) \mid i=0,\ldots,q-2\}.
\end{align*}
Their union is a difference set of the form~\eqref{twinds2}, because $-1$ is a square in
one of the fields $\F_q$, $\F_{q+2}$ and a non-square in the other. By subtracting two non-zero
coordinates $x$ and $y$, we get
\begin{align*}
 \{(\alpha^{2i+x}(1-\alpha^{y-x}),\beta^{2j+x}(1-\beta^{y-x})) \mid & \,\, i=0,\ldots,(q-3)/2,\\
                                                                    & \,\, j=0,\ldots,(q-1)/2\},\\[1mm]
\{(\alpha^{2i+x+1}(1-\alpha^{y-x}),\beta^{2j+x+1}(1-\beta^{y-x})) \mid & \,\, i=0,\ldots,(q-3)/2,\\
                                                                          & \,\, j=0,\ldots,(q-1)/2\},\\[1mm]
 \{(0,0)\}\cup\{(\alpha^{i+y}(1-\alpha^{y-x}),0) \mid & \,\, i=0,\ldots,q-2\}.
\end{align*}
The union is a difference set of the form~\eqref{twinds1} if the elements $1-\alpha^{y-x}$ and
$1-\beta^{y-x}$ are both squares or both non-squares. If one is a square, and
the other a non-square, we get a difference set of the form~\eqref{twinds2}. Hence, $D$ is a
$q$-dimensional $(4m-1,2m-1,m-1)$ difference set.
\end{proof}

The constructions of Theorems~\ref{tmpaley} and~\ref{tmtpp} are implemented as\break PAG~\cite{PAG}
commands \texttt{PaleyDifferenceSet} and \texttt{Twin\-Prime\-Power\-Dif\-fe\-rence\-Set}. It would be of
interest to extend other known families of difference sets to higher dimensions.

\section{Examples with small parameters}\label{sec5}

Let $\mu_G(v,k,\lambda)$ be the largest integer~$n$ such that an $n$-dimensional
$(v,k,\lambda)$ difference set in~$G$ exists. We write $\mu(v,k,\lambda)$ if
there is only one group of order~$v$ up to isomorphism. Proposition~\ref{dzprop}
implies $\mu_G(v,k,\lambda)\le \nu(v,k,\lambda)$. Theorems from Section~\ref{sec4}
give lower bounds on~$\mu$ for particular parameters and elementary abelian
groups, e.g.\ Theorem~\ref{tmpaley} implies $\mu_G(q,(q-1)/2,(q-3)/4)\ge q$ for
$q=p^s$ and $G=(\Z_p)^s$.

For small parameters $(v,k,\lambda)$, we can determine the
exact value of~$\mu$ by performing a complete classification of
$n$-dimensional difference sets. For example, each of the three
$2$-subsets $\{0,1\}$, $\{0,2\}$, $\{1,2\}$ is a $(3,2,1)$
difference set in~$\Z_3$. We add a zero coordinate to get
$2$-dimensional difference sets $\{(0,0),(0,1)\}$, $\{(0,0),(0,2)\}$,
$\{(0,1),(0,2)\}$ and try to extend them by adding a third
coordinate. Our difference sets are normalized, hence the added coordinate
must contain a permutation of a difference set in~$\Z_3$. There are
six possible extensions for each $2$-dimensional difference set,
but only three of them satisfy the condition from Definition~\ref{defdifset}.
We get a total of nine $3$-dimensional difference sets. Their
developments are equivalent, so we can continue extending only
one of them, e.g.\ $\{(0,0,1),(0,1,0)\}$. Now none of the six
extensions satisfy the condition for a $4$-dimensional difference
set, and we conclude that $\mu(3,2,1)=3$.

This procedure was performed in GAP~\cite{GAP}, using the catalog of
difference sets from the DifSets package~\cite{DP19}. The results are
summarized in Table~\ref{tab2}. There is no simple way of complementing
higher-dimensional difference sets, so the calculation was performed
for each of the complementary parameters $(v,k,\lambda)$ and
$(v,v-k,v-2k+\lambda)$, when possible. The values of $\mu(15,7,3)$ and
$\mu(15,8,4)$ are indeed different. For each increase of dimension
by one, the number of candidates to consider is proportional to the
number of $(v,k,\lambda)$ difference sets times~$k!$. Because of the
factor~$k!$, we could perform a complete classification for $(13,4,1)$
and $(21,5,1)$, but not for the complementary parameters $(13,9,6)$
and $(21,16,12)$. The result $\mu_G(21,5,1)=3$ is valid for both groups
$G=\Z_{21}$ and $G=\Z_7 \rtimes \Z_3$.

\begin{table}[!h]
\caption{Maximal dimensions of small difference sets.}\label{tab2}
\begin{tabular}{cccc}
\toprule
$\mu(7,3,1)=7$ & $\mu(11,5,2)=11$ & $\mu(15,7,3)=3$ & $\mu(13,4,1)=13$ \\[1mm]
$\mu(7,4,2)=7$ & $\mu(11,6,3)=11$ & $\mu(15,8,4)=4$ & $\mu_G(21,5,1)=3$ \\
\bottomrule
\end{tabular}
\end{table}

We were also able to classify higher-dimensional $(16,6,2)$ difference
sets in some groups. Twelve of the $14$ groups of order~$16$
contain difference sets~\cite{DP19}. We computed exact values of~$\mu_G(16,6,2)$
for~$8$ of these groups, and lower bounds for~$4$ groups. The results are
summarized in Table~\ref{tab3}. The second column contains a description
of the structure of~$G$. This does not identify the group uniquely, so in
the first column we give the ID of~$G$ in the GAP~\cite{GAP} library of
small groups. The column Nds contains the number of difference
sets in~$G$ up to equivalence, and Tds contains the total number of difference
sets. The last column contains the value of $\mu_G(16,6,2)$, or a lower bound.

\begin{table}[t]
\caption{Maximal dimensions of $(16,6,2)$ difference sets.}\label{tab3}
\begin{tabular}{ccccc}
\toprule
ID & $G$ & Nds & Tds & $\mu_G(16,6,2)$\\
\midrule
1 & $\Z_{16}$ & 0 & 0 & -- \\
2 & $\Z_4^2$ & 3 & 192 & 4 \\
3 & $(\Z_4 \times \Z_2) \rtimes \Z_2$ & 4 & 192 & 4 \\
4 & $\Z_4 \rtimes \Z_4$ & 3 & 192 & 4 \\
5 & $\Z_8 \times \Z_2$ & 2 & 192 & 4 \\
6 & $\Z_8 \rtimes \Z_2$ & 2 & 64 & 4 \\
7 & $D_{16}$ & 0 & 0 & -- \\
8 & $QD_{16}$ & 2 & 128 & 4 \\
9 & $Q_{16}$ & 2 & 256 & 4 \\
10 & $\Z_4 \times \Z_2^2$ & 2 & 448 & $\ge 4$\\
11 & $\Z_2 \times D_8$ & 2 & 192 & 4 \\
12 & $\Z_2 \times Q_8$ & 2 & 704 & $\ge 8$ \\
13 & $(\Z_4 \times \Z_2) \rtimes \Z_2$ & 2 & 320 & $\ge 6$ \\
14 & $\Z_2^4$ & 1 & 448 & $\ge 14$ \\
\bottomrule
\end{tabular}
\end{table}

The classification of $n$-dimensional difference sets produced a lot of
examples of projection $n$-cubes. Lower bounds on the number of inequivalent
cubes in $\P^n(v,k,\lambda)$ are given in Table~\ref{tab4}. In most cases,
these are all possible examples obtained from difference sets, except
for $(16,6,2)$ and $n\ge 4$. An online version of the
table with links to files containing the examples in GAP format
is available at
\begin{center}
\url{https://web.math.pmf.unizg.hr/~krcko/results/pcubes.html}
\end{center}

\begin{table}[!b]
\caption{Lower bounds on the number of $\P^n(v,k,\lambda)$-cubes.}\label{tab4}
{\small \begin{tabular}{c|cccccccccccc}
\toprule
 & \multicolumn{12}{c}{$n$} \\[1mm]
$(v,k,\lambda)$ & 3 & 4 & 5 & 6 & 7 & 8 & 9 & 10 & 11 & 12 & 13 & 14 \\
\midrule
$(7,3,1)$ & 2 & 2 & 1 & 1 & 1 & & & & & & & \\[1mm]
$(7,4,2)$ & 2 & 2 & 1 & 1 & 1 & & & & & & & \\[1mm]
$(11,5,2)$ & 2 & 4 & 6 & 6 & 4 & 2 & 1 & 1 & 1 &  &  &  \\[1mm]
$(11,6,3)$ & 2 & 4 & 6 & 6 & 4 & 2 & 1 & 1 & 1 &  &  &  \\[1mm]
$(13,4,1)$ & 3 & 7 & 10 & 14 & 14 & 10 & 7 & 3 & 1 & 1 & 1 &  \\[1mm]
$(15,7,3)$ & 3 &  &  &  &  &  &  &  &  &  &  &  \\[1mm]
$(15,8,4)$ & 6 & 1 &  &  &  &  &  &  &  &  &  &  \\[1mm]
$(16,6,2)$ & 724 & 8464 & 1601 & 1754 & 986 & 505 & 178 & 70 & 16 & 7 & 2 & 1 \\[1mm]
$(21,5,1)$ & 6 &  &  &  &  &  &  &  &  &  &  &  \\
\bottomrule
\end{tabular}}
\end{table}

So far, the only examples of projection cubes not coming from
difference sets are $C_5\in \P^5(3,2,1)$ from Example~\ref{exC45}
and its restrictions to dimensions~$3$ and~$4$. The $3$-cube~$C_4$
is equivalent to the development of $\{(0,0,1),(0,1,0)\}$
over~$\Z_3$. We found more examples not coming from difference
sets in $\P^3(16,6,2)$.

Firstly, we computed the full autotopy groups of the $724$ cubes
$C\in \P^3(16,6,2)$ obtained from difference sets.
We then chose subgroups $G\le \Atop(C)$ of orders
$|G|\in \{18, 24, 32, 36, 48, 60\}$ and constructed all
$\P^3(16,6,2)$-cubes with $G$ as a prescribed autotopy group.
We used a modification of the Kramer-Mesner method~\cite{KM76},
similar as in~\cite{KPT24}. For each choice of~$G\le \Atop(C)$, we always
get the cube~$C$, but sometimes we also find other inequivalent cubes.
In this way, we constructed $102$ new examples in $\P^3(16,6,2)$. They
are not equivalent to cubes coming from difference sets because their
full autotopy groups act non-transitively on the coordinates
(cf.\ Proposition~\ref{difatop}). The orbit size distributions that occurred
are $1+6+9$, $4+12$, $6+10$, and $8+8$. An improved lower bound for
the number of inequivalent $\P^3(16,6,2)$-cubes is $724+102=826$.

\begin{figure}[!b]
\begin{center}
\includegraphics[width=60mm]{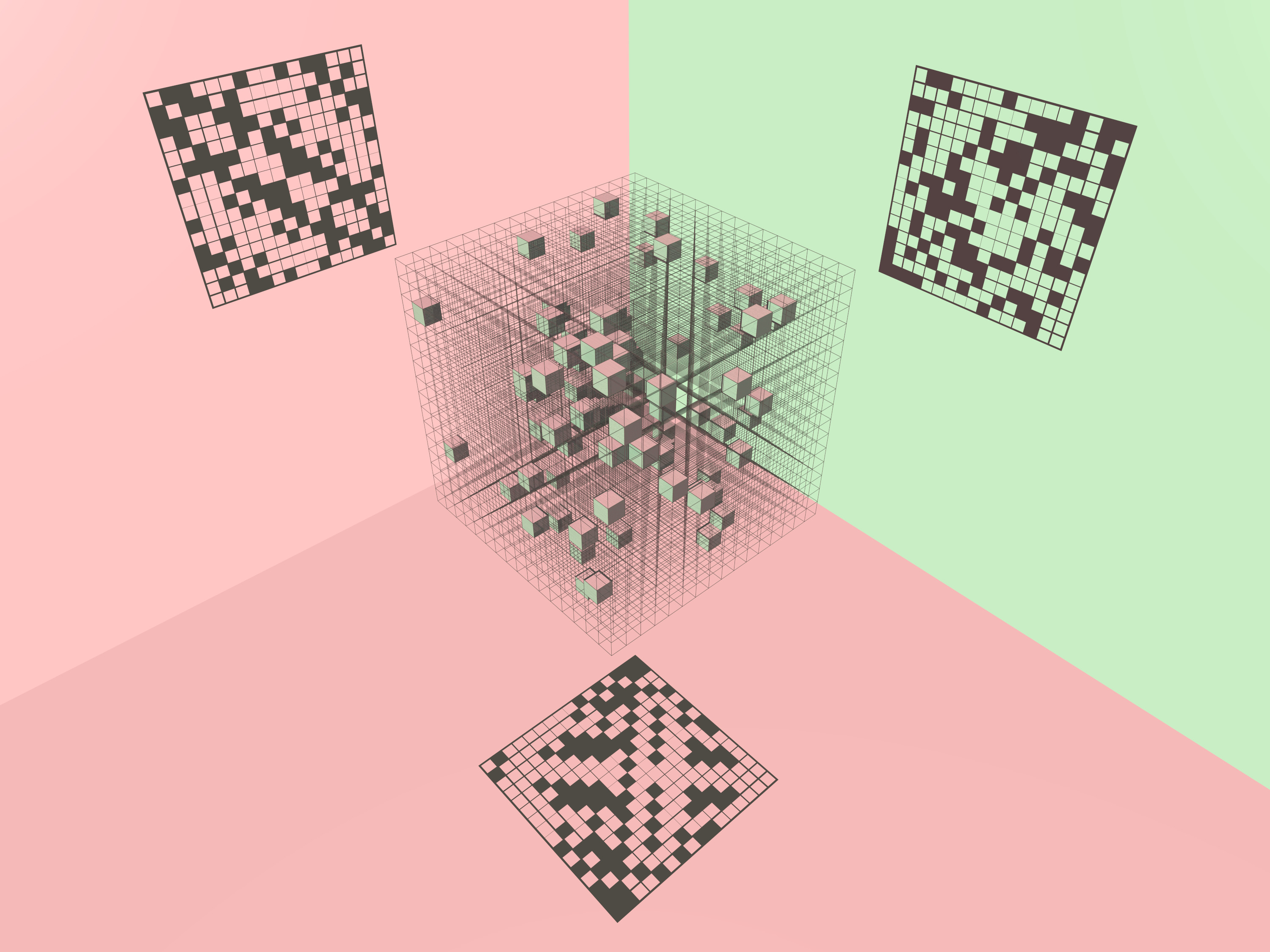}\hskip 4mm\includegraphics[width=60mm]{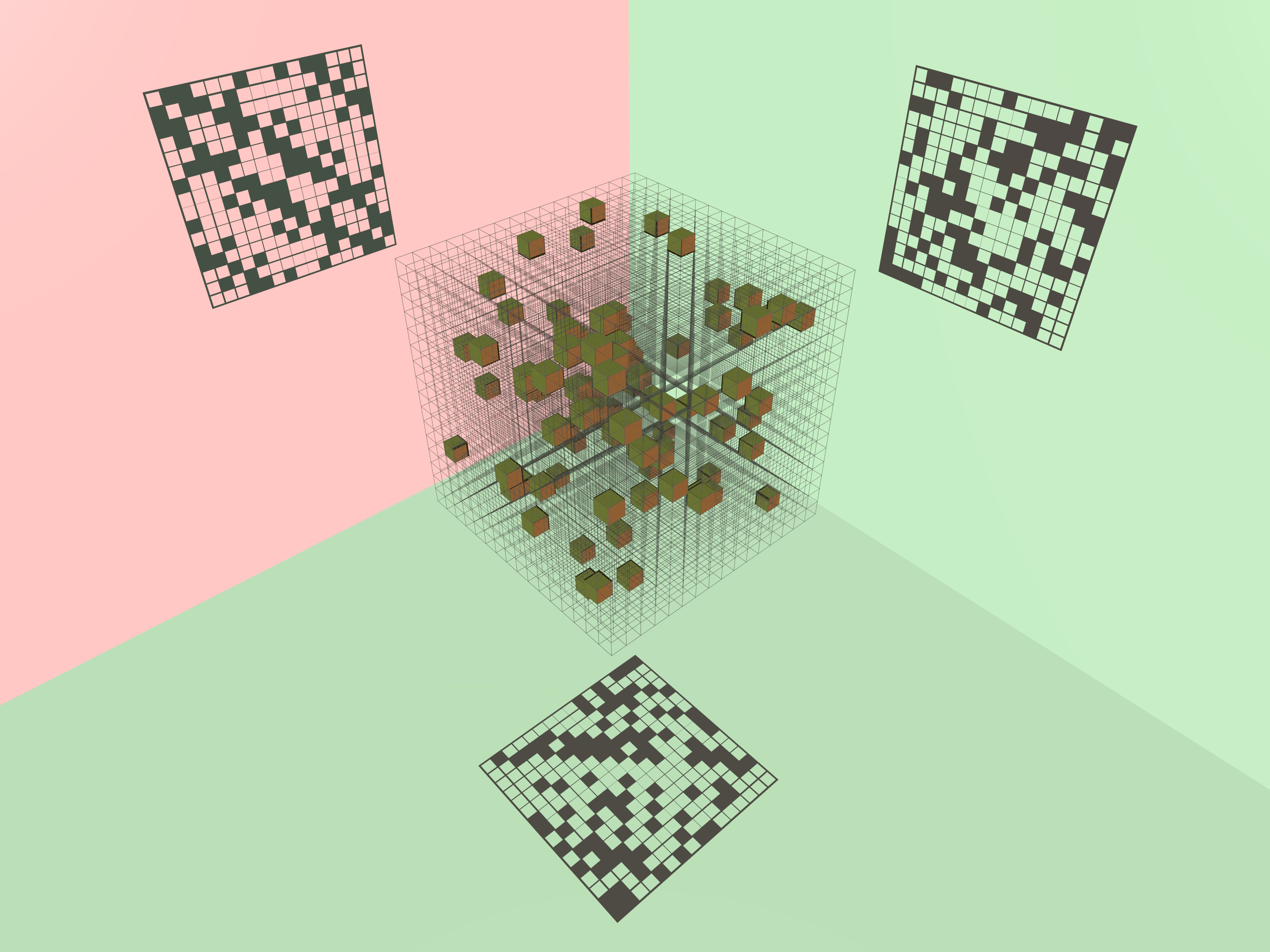}\\[4mm]
\includegraphics[width=60mm]{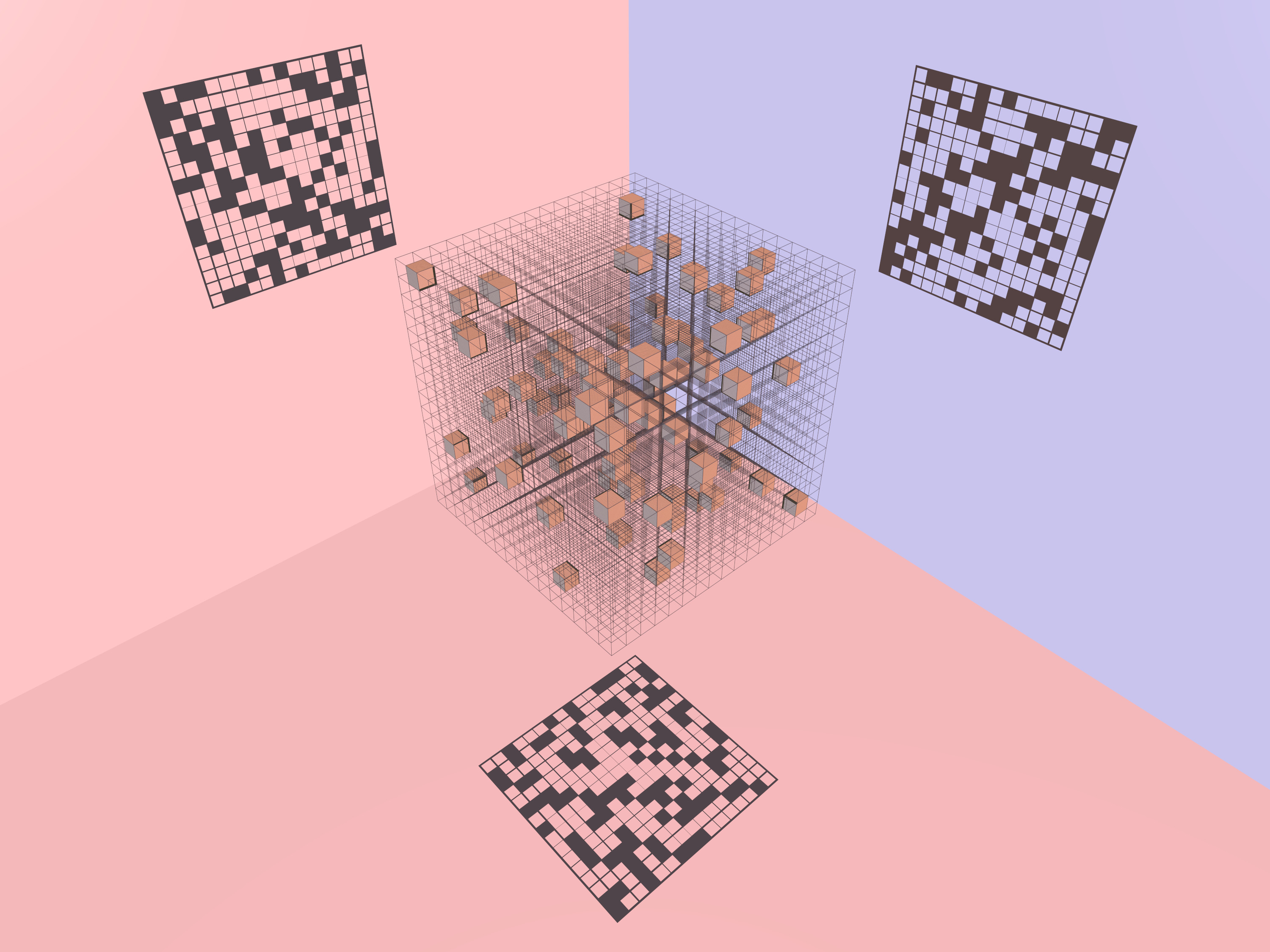}\hskip 4mm\includegraphics[width=60mm]{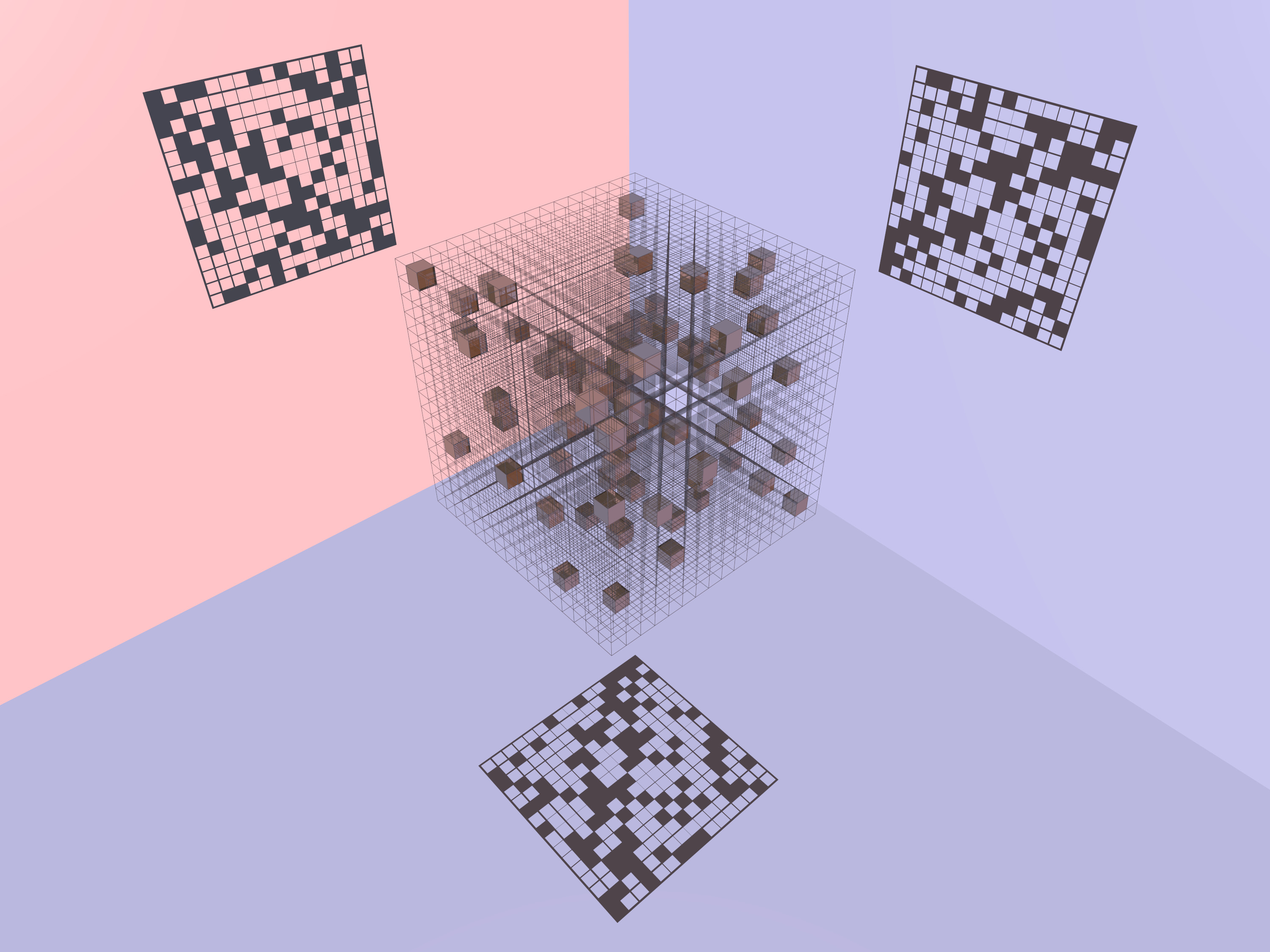}
\end{center}
\caption{Cubes in $\P^3(16,6,2)$ with non-isomorphic projections.}\label{fig2}
\end{figure}

The $102$ non-difference cubes are available on our
\href{https://web.math.pmf.unizg.hr/~krcko/results/pcubes.html}{web page}. Four
of them are shown in Figure~\ref{fig2} and their GAP representations
are available in the file
\href{https://web.math.pmf.unizg.hr/~krcko/results/pcubes/examples.oa}{\texttt{examples.oa}}.
These four cubes are interesting because they have non-isomorphic projections.
Up to isomorphism, there are exactly three $(16,6,2)$ designs~\cite{MR07}.
We denote them $\D_R$, $\D_G$, $\D_B$ and call them the \emph{red}, \emph{green},
and \emph{blue} design, respectively. They are identified by orders
of the full automorphism groups: $|\Aut(\D_R)|=11520$,
$|\Aut(\D_G)|=768$, and $|\Aut(\D_B)|=384$. Figure~\ref{fig2} was rendered
with red, green, and blue light sources so that the shadows appear in
the appropriate color. Thus, the four cubes have projections
$(\D_R,\D_R,\D_G)$, $(\D_R,\D_G,\D_G)$, $(\D_R,\D_R,\D_B)$, and
$(\D_R,\D_B,\D_B)$.

An interesting question is  whether there exists a projection cube with
all three designs $\D_R$, $\D_G$, and $\D_B$ appearing as projections. So far
we have not been able to find such an example.

\section{Concluding remarks}\label{sec6}

We have introduced a new type of generalization of symmetric block designs
to higher dimensions, called \emph{projection cubes}. The defining property
is that all $2$-dimensional projections are symmetric $(v,k,\lambda)$ designs.
This leads to bounds on the dimension~$n$ in terms of~$v$ and~$k$.
A previously studied generalization was based on $2$-dimensional sections,
and the dimension was not bounded for given parameters $(v,k,\lambda)$.
We have also defined $n$-dimensional difference sets suitable for
constructing projection cubes, and extended some classic results about
difference sets to higher dimensions. Infinite families of these objects
are constructed in Theorems~\ref{tmpaley} and~\ref{tmtpp}. Some sporadic
examples are constructed by calculations in GAP~\cite{GAP}.

We are currently working on an implementation of the classification algorithm
for $n$-dimensional difference sets in the C programming language. We want to
extend the computational results presented in Tables~\ref{tab2}--\ref{tab4} to
larger parameters $(v,k,\lambda)$. Future research directions include generalizing
other known results about symmetric designs and difference sets to higher dimensions.
In particular, it would be of interest to find more examples of projection cubes
not associated with difference sets. For some parameters such as $(25,9,3)$, symmetric
designs exist~\cite{MR07}, but there are no difference sets~\cite{DP19}.
The existence of projection $n$-cubes with these parameters is open for $n\ge 3$.

\vskip 2mm

\textsc{Acknowledgements.}

The authors wish to thank an anonymous referee who pointed out the improved
bound of Theorem~\ref{dimbound2}. We also thank the editors for suggestions
improving the presentation of the paper.


\begin{thebibliography}{20}

\bibitem{LB71}
Baumert, L.D.: Cyclic difference sets.
Lecture Notes in Math.\ 182. Springer-Verlag,
Berlin, New York (1971). \url{https://doi.org/10.1007/BFb0061260}

\bibitem{BJL99}
Beth, T., Jungnickel, D., Lenz, H.: Design theory. Second
edition. Volume~I. Cambridge University Press, Cambridge (1999).
\url{https://doi.org/10.1017/CBO9780511549533}

\bibitem{RHB55}
Bruck, R.H.: Difference sets in a finite group,
Trans.\ Amer.\ Math.\ Soc. 78, 464--481 (1955).
\url{https://doi.org/10.1090/S0002-9947-1955-0069791-3}

\bibitem{dL90}
de~Launey, W.: On the construction of $n$-dimensional designs from $2$-dimensional
designs. Australas.\ J.\ Combin. 1, 67--81 (1990).
\url{https://ajc.maths.uq.edu.au/pdf/1/ocr-ajc-v1-p67.pdf}

\bibitem{dLH93}
de Launey, W., Horadam, K.J.: A weak difference set construction for
higher-dimensional designs. Des.\ Codes Cryptogr.\ 3, 75--87 (1993).
\url{https://doi.org/10.1007/BF01389357}

\bibitem{JHD07}
Dinitz, J.H.: \emph{Room squares}. In: Colbourn, C.J., Dinitz, J.H. (eds.)
Handbook of combinatorial designs. Second edition, pp.\ 584--590.
Chapman \& Hall/CRC, Boca Raton, FL (2007).
\url{https://doi.org/10.1201/9781420010541}

\bibitem{GAP}
The GAP Group: GAP -- Groups, Algorithms, and Programming, Version 4.13.1
(2024). \url{https://www.gap-system.org}

\bibitem{GC07}
Greig, M., Colbourn, C.J.: Orthogonal arrays of index more than one.
In: Colbourn, C.J., Dinitz, J.H. (eds.) Handbook of combinatorial designs.
Second edition, pp.\ 219--224. Chapman \& Hall/CRC, Boca Raton, FL (2007).
\url{https://doi.org/10.1201/9781420010541}

\bibitem{HSS99}
Hedayat, A.S., Sloane, N.J.A., Stufken, J.: Orthogonal arrays.
Theory and applications. Springer-Verlag, New York (1999).
\url{https://doi.org/10.1007/978-1-4612-1478-6}

\bibitem{IT07}
Ionin, Y.J., Tran, V.T.: \emph{Symmetric designs}. In: Colbourn, C.J., Dinitz, J.H. (eds.)
Handbook of combinatorial designs. Second edition, pp.\ 110--124. Chapman \& Hall/CRC,
Boca Raton, FL (2007). \url{https://doi.org/10.1201/9781420010541}

\bibitem{JPS07}
Jungnickel, D., Pott, A., Smith, K.W.: Difference sets. In: Colbourn, C.J.,
Dinitz, J.H. (eds.) Handbook of combinatorial designs. Second edition,
pp.\ 419--435. Chapman \& Hall/CRC, Boca Raton, FL (2007).
\url{https://doi.org/10.1201/9781420010541}

\bibitem{DEK92}
Knuth, D.E.: Two notes on notation. Amer.\ Math.\ Monthly 99, 403--422 (1992).
\url{https://doi.org/10.1080/00029890.1992.11995869}

\bibitem{KM76}
Kramer, E.S., Mesner, D.M.: $t$-designs on hypergraphs.
Discrete Math.\ 15, 263--296 (1976).
\url{https://doi.org/10.1016/0012-365X(76)90030-3}

\bibitem{PAG}
Kr\v{c}adinac, V.: PAG -- Prescribed Automorphism Groups, Version 0.2.4
(2024). \url{https://doi.org/10.5281/zenodo.14194448}

\bibitem{KPT24}
Kr\v{c}adinac, V., Pav\v{c}evi\'{c}, M.O., Tabak, K.: Cubes of symmetric
designs. Ars Math.\ Contemp.\ 25, \#P1.10 (2025).
\url{https://doi.org/10.26493/1855-3974.3222.e53}

\bibitem{MR07}
Mathon, R., Rosa, A.: $2$-$(v,k,\lambda)$ designs of small order. In: Colbourn, C.J.,
Dinitz, J.H. (eds.) Handbook of combinatorial designs. Second edition,
pp.\ 25--58. Chapman \& Hall/CRC, Boca Raton, FL (2007).
\url{https://doi.org/10.1201/9781420010541}

\bibitem{DP19}
Peifer, D.: \emph{DifSets, an algorithm for enumerating all difference
sets in a group}, Version 2.3.1 (2019). \url{https://dylanpeifer.github.io/difsets}

\bibitem{POVRay}
Persistence of Vision Raytracer, Version 3.7 (2013). Persistence of Vision Pty.
Ltd., Williamstown, Victoria, Australia. \url{http://www.povray.org/}

\end{thebibliography}
\end{document}